\documentclass[a4paper,12pt,reqno]{amsart}
\usepackage{amssymb}
\usepackage[dvips]{graphicx}
\usepackage{hyperref}

\usepackage[dvipsnames]{xcolor} 
\usepackage{colortbl}
\usepackage{tikz-cd}

\usepackage{bigdelim}
\usepackage{arydshln}

\input epsf.tex
\newdimen\xsize
\newdimen\oldbaselineskip
\newdimen\oldlineskiplimit
\def\restorelineskip{\baselineskip=\oldbaselineskip%
\lineskiplimit=\oldlineskiplimit}
\def\putm[#1][#2]#3{
\hbox{\vbox to 0pt{\parindent=0pt%
\vskip#2\xsize\hbox to0pt{\hskip#1\xsize $#3$\hss}\vss}}}%
\long\def\Line#1{\hbox to \hsize{#1}}
\def\putt[#1][#2]#3{
\vbox to 0pt{\noindent\hskip#1\xsize\lower#2\xsize%
\vtop{\restorelineskip#3}\vss}}

\makeatletter
\def\xbig[#1]#2{{\hbox{$\m@th\left#2\vbox to#1\xsize{}%
\right.\n@space$}}}
\makeatother
\def\xlar[#1]#2{%
\smash{\mathop{ \hbox to #1\xsize{\leftarrowfill}}\limits^{#2}}}
\def\xrar[#1]#2{%
\smash{\mathop{ \hbox to #1\xsize{\rightarrowfill}}\limits^{#2}}}
\def\xline[#1]{\hbox to #1\xsize{\leaders\hrule\hfill}}

\DeclareFontFamily{U}{rsf}{\skewchar\font'177}%
\DeclareFontShape{U}{rsf}{m}{n}{<-6>rsfs5<6-8>rsfs7<8->rsfs10}{}%
\DeclareFontShape{U}{rsf}{b}{n}{<-6>rsfs5<6-8>rsfs7<8->rsfs10}{}%
\DeclareMathAlphabet\RSFS{U}{rsf}{m}{n}
\SetMathAlphabet\RSFS{bold}{U}{rsf}{b}{n}
\let\scr=\rfs

\def\sf#1{{\mathsf{#1}}}

\def\slsf{\slshape \sffamily }

\def\msmall#1{\mathchoice{\hbox{\small$\displaystyle {#1}$}}{#1}{#1}{#1}}

\def\cc{{\mathbb C}}

\def\nn{{\mathbb N}}

\def\pp{{\mathbb P}}

\def\deg{\sf{deg}\,}

\def\lim{\mathop{\sf{lim}}}

\def\mto{\dashrightarrow}

\def\pr{\sf{pr}}

\def\<{\langle}\let\la=\<
\def\>{\rangle}\let\ra=\>

\def\ddef{\mathrel{{=}\raise0.3pt\hbox{:}}}
\def\deff{\mathrel{\raise0.3pt\hbox{\rm:}{=}}}

\def\fraction#1/#2{\mathchoice{{\msmall{ #1\over#2}}}%
{{ #1\over #2 }}{{#1/#2}}{{#1/#2}}}
\def\norm#1{\left\Vert{#1}\right\Vert}

\def\le{\leqslant}

\def\emptyset{\varnothing}

\def\longpoints{\leaders\hbox to 0.5em{\hss.\hss}\hfill \hskip0pt}
\def\stateskip{\smallskip}
\def\state#1. {\stateskip\noindent{\bf#1. }} 
\def\statep#1. {\stateskip\noindent{\bf#1 }} 
\def\proof{\state Proof. \2}

\def\Chi{\raise 2pt\hbox{$\chi$}}

\def\ie{\hskip1pt plus1pt{\sl i.e.\/,\ \hskip1pt plus1pt}}

\def\Chi{\raise 2pt\hbox{$\chi$}}

\let\phI=\phi\let\phi=\varphi\let\varphi=\phI
\def\scra{\scr{A}}

\let\cal=\mathcal
\def\cala{{\cal A}}
\def\calb{{\cal B}}

\def\cald{{\cal D}}

\def\calh{{\cal H}}
\def\cali{{\cal I}}
\def\calj{{\cal J}}

\def\calm{{\cal M}}

\def\calo{{\cal O}}
\def\calp{{\cal P}}

\def\calr{{\cal R}}

\def\calu{{\cal U}}
\def\calv{{\cal V}}

\def\calx{{\cal X}}
\def\caly{{\cal Y}}

\def\1{{1\mkern-5mu{\rom l}}}

\def\ge{\geqslant}

\def\fraction#1/#2{\mathchoice{{\msmall{ #1\over#2}}}%
{{ #1\over #2 }}{{#1/#2}}{{#1/#2}}}

\def\le{\leqslant}

\def\emptyset{\varnothing}
\newcommand{\2}{\thinspace}

\def\qed{\ \ \hfill\hbox to .1pt{}\hfill\hbox to .1pt{}\hfill $\square$\par}

\def\comment#1\endcomment{}

 \belowdisplayskip=\abovedisplayskip

\def\lineeqqno(#1){\hfill\llap{\vbox to 10pt%
{\vss\begin{align} \eqqno(#1)\end{align}\vss}}\vskip1pt}

\advance\headheight 1.2pt

\def\ShowwLLabel#1{}

\def\thechpt{\Roman{chpt}}

\def\newchapt[#1]#2{\newpage%
\refstepcounter{chpt}\setcounter{subsection}{0}%
\setcounter{thm}{0}\setcounter{defi}{0}%
\setcounter{rema}{0}\setcounter{exrc}{0}%
\renewcommand{\thesubsection}{\thechpt.\arabic{subsection}}%
\section*{\begin{center}\huge \bf Chapter \thechpt\\
#2 \end{center}}\label{#1}%
\ \smallskip%
\markboth{Chapter \thechpt}{#2}%
}
\def\newsect[#1]#2{\refstepcounter{section}\setcounter{equation}{0}%
\renewcommand{\thesubsection}{\arabic{section}.\arabic{subsection}}%
\section*{\arabic{section}.
#2}\vspace{-20pt}\label{#1}\vspace{20pt}%
\markboth{Section \arabic{section}}{#2}}

\def\newlect[#1]#2{\refstepcounter{section}%
\renewcommand{\thesubsection}{\arabic{section}.\arabic{subsection}}%
\section*{Lecture \arabic{section}\\
#2}\label{#1}%
\markboth{Lecture \arabic{section}}{#2}}

\def\newprg[#1]#2{\refstepcounter{subsection}%
\subsection*{{\thesubsection.\ #2}} \label{#1}%
}

\def\newappx[#1]#2{%
\refstepcounter{appx}\setcounter{section}{0}%
\renewcommand{\thesubsection}{A\arabic{appx}.\arabic{subsection}}%
\section*{Appendix \arabic{appx}\\ #2}
\label{#1}%
\markboth{Appendix A\arabic{appx}}{#2}
}

\newtheorem{mainthm}{Main Theorem}
   \def\newthm#1{\begin{mainthm}}
   \def\refthm#1{{\sl Main Theorem }}

\newtheorem{thm}{Theorem}[section]
   \def\newthm#1{\begin{thm}\label{#1}}
   \def\refthm#1{{\sl Theorem \ref{#1}}}

\newtheorem{nnthm}{Theorem}
   \def\newthm#1{\begin{nnthm}\label{#1}}
   \def\refthm#1{{\sl Theorem \ref{#1}}}
\newtheorem{lem}{Lemma}[section]
   \def\newlemma#1{\begin{lem} \label{#1}}

\newtheorem{prop}{Proposition}[section]
   \def\newprop#1{\begin{prop}\label{#1}}

\newtheorem{nnprop}{Proposition}
   \def\newprop#1{\begin{nnprop}\label{#1}}

\newtheorem{corol}{Corollary}[section]
   \def\newcorol#1{\begin{corol} \label{#1}}

\newtheorem{nncorol}{Corollary}
   \def\newcorol#1{\begin{nncorol} \label{#1}}

\newtheorem{defi}{Definition}[section]
   \def\newdefi#1{\begin{defi} \label{#1}\rm }

\newtheorem{nndefi}{Definition}
   \def\newdefi#1{\begin{nndefi} \label{#1}\rm }

\newtheorem{exmp}{Example}[section]
   \def\newexmp#1{\begin{exmp} \label{#1}\rm }

\newtheorem{nnexmp}{Example}
   \def\newthm#1{\begin{nnexmp}\label{#1}}
   \def\refthm#1{{\sl Example \ref{#1}}}
\newtheorem{exrc}{Exercise}
   \def\newexrc#1{\begin{exrc} \label{#1}\rm }

\newtheorem{rema}{Remark}[section]
   \def\newrema#1{\begin{rema} \label{#1}\rm }

\newtheorem{nnrema}{Remark}
   \def\newrema#1{\begin{nnrema} \label{#1}\rm }

\def\eqqno(#1){\label{(#1)}}
\def\eqqref(#1){(\ref{(#1)})}

\def\el2{\sf{L^2}}
\def\l{\it{l}}
\def\vect{{\mathrm{v}}}
\def\wect{{\mathrm{w}}}

\title{A note on meromorphicity in infinite dimensions}
\author{Ivashkovich S.}
\address{Universit\'e de Lille, UFR de Math\'ematiques, 59655 Villeneuve
d'Ascq, France.}
\email{serge.ivashkovych@univ-lille.fr}

\thanks{Partially supported by the Labex Cempi ANR-11-LABX-0007-01}
\date{\today}

\begin{document}

\begin{abstract}
The goal of this note is to develop the notion of a meromorphic mapping 
between infinite-dimensional complex manifolds. We restrict ourselves with 
second countable  Hilbert manifolds modelled on a separable Hilbert space.  
\end{abstract}
\maketitle

\setcounter{tocdepth}{1}
\tableofcontents

\newsect[INT]{Introduction and statement of the main result}

\newprg[INT.gen-def]{General definition}
We begin by recalling the notion of a meromorphic map between finite 
dimensional complex manifolds. Let $X$ and $Y$ be such manifolds.

\begin{nndefi}
\label{fin-dim}
A meromorphic map $\hat f: X\mto Y$ from $X$ to $Y$
if defined by an {\slsf irreducible} analytic set $\Gamma_{\hat f}$ 
in $X\times Y$, called the graph of $\hat f$, such that:

\begin{enumerate}
\item[1)]the restriction to  $\Gamma_{\hat f}$ of the natural projection $\pi_X:X\times Y\to X$  is proper;

\smallskip\item[2)] there exists an analytic subset $I$ of $X$ of codimension 
$\ge 2$ such that over $X\setminus I$ the set $\Gamma_{\hat f}$ is the 
graph of a holomorphic map.
\end{enumerate}
\end{nndefi}

\smallskip In other words there exists a holomorphic map $f:X\setminus I\to Y$
such that $\Gamma_{\hat f}\cap \left[\left(X\setminus I\right)\times Y\right]$ 
is the graph $\Gamma_f$ of $f$. The minimal $I$ satisfying condition (2)
of this definition is called the {\slsf indeterminacy set} of $\hat f$ and 
will be denoted as $I_{\hat f}$. It is worth noting that when $Y=\pp^n$, a
complex projective space, the data of such $\Gamma_{\hat f}$ is 
equivalent to the following:

\smallskip {\it  
\begin{enumerate}
\item[1)] given meromorphic functions $f_1,...,f_n$ on $X$, set
$I_{\hat f} = \bigcup_{k=1}^n I_{f_k}$, where $I_{f_k}$ stands for the 
indeterminacy set of $f_k$;

\smallskip\item[2)]  for every point $a\in X$ there exists a neighbourhood
$V\ni a$ and holomorphic in $V$ functions $\phi_0,...,\phi_n$ without common factors such that $f_k=\frac{\phi_k}{\phi_0}$ on $V$, $k=1,...,n$;

\smallskip\item[3)]  on $V\setminus I_{\hat f}$ one has $f(z) = [\phi_0(z):...:
\phi_n(z)]$ and, moreover, $ I_{\hat f}\cap V = \{z\in V:\phi_0(z) = ...=\phi_n(z)=0\}$.
\end{enumerate}
}

\smallskip Our task in this note is to propose  a definition of a meromorphic
map between complex Hilbert manifolds. We restrict ourselves to the case of 
complex Hilbert manifolds that are modelled on a separable complex Hilbert
space $L$, \ie on $\cc^n$ or $\l^2$. In addition throughout this paper, if the opposite is not explicitly stated,
we suppose our manifolds to be second countable, i.e. their topology has a countable base. Notice that second countable Hilbert manifolds are paracompact, admit a partition of unity and therefore a Hermitian metric. First of all let us 
point out that Definition \ref{fin-dim} doesn't give the right notion 
in infinite dimension.

\begin{nnexmp}
\label{hilb-proj}
\smallskip\noindent{\slsf Projection $\l^2\setminus\{0\}\to \pp (\l^2)$.} \rm
Recall that the Hilbert projective space is $\pp (\l^2)\deff (l^2\setminus\{0\})/w
\sim \lambda w$, $\lambda \not= 0$, where $w=(w_0,w_1,...)\in \l^2$. By $z^j_k =
w_j/w_k$ we denote the affine coordinates in the affine chart $\calu_k = 
\{[w]:w_k\not=0\} = \pp (l^2)\setminus \calh_k$, where $\calh_k$ is the hyperplane 
\begin{equation}
\eqqno(hyp-hj)
\calh_k\deff \{[w]\in \pp (\l^2): w_k=0\}. 
\end{equation}
The graph $\Gamma_{\pr}\subset \left(\l^2\setminus \{0\}\right)\times \pp (\l^2)$
of the natural projection 
\[
\pr : \l^2\setminus \{0\}\to \pp (\l^2),
\]
\begin{equation}
\eqqno(bl-up1)
\l^2\setminus \{0\}\ni (z_0,.z_1,...) = z \to [z] \deff [z_0:z_1:...] \in 
\pp (\l^2)
\end{equation}
closes in $\l^2\times \pp (\l^2)$ to the following analytic set:
\begin{equation}
\eqqno(graph-proj)
\Gamma_{\widehat{\pr}} = \{(z,[w])\in \l^2\times \pp (\l^2): w_1z_0-w_0z_1
= ...= w_kz_0-w_0z_k = ... =0\},
\end{equation}
and this certainly should be the graph of a meromorphic map.
But notice that the preimage of the point $0$ under (another) natural projection 
$\pi :\Gamma_{\widehat{\pr}}\to \l^2$  is the whole $\pp (\l^2)$, 
\ie is {\slsf not compact!}
\end{nnexmp}

Let $\calx$ and $\caly$ be hermitian Hilbert manifolds, \ie complex Hilbert 
manifolds equipped with a Hermitian metric. We suggest to introduce the 
following notion.

\begin{nndefi}
\label{inf-dim}
A meromorphic mapping $\hat f: \calx \mto \caly$ is defined by an 
irreducible analytic set $\Gamma_{\hat f}$ in $\calx\times \caly$ such that:
\begin{enumerate}
\item[1)] for every bounded and complete subset $B \subset \calx$ the set 
$\pi_{\calx}^{-1}(B)\cap \Gamma_{\hat f}$ is bounded and complete in 
$\Gamma_{\hat f}$;
\item[2)] there exists an analytic set $\cali$ in $\calx$ of codimension 
$\geqslant 2$ and a holomorphic map $f: \calx \setminus \cali \to \caly$
such that $\Gamma_{\hat f}\cap \big[\left(\calx\setminus \cali\right)\times \caly\big]$ is the graph $\Gamma_f$ of $f$.
\end{enumerate} 
\end{nndefi}

Here $\pi_{\calx}: \calx \times \caly \to \calx$ as above stands for the 
natural 
projection. The minimal $\cali$ such that $\hat f$ is holomorphic  on 
$\calx\setminus \cali$ is called the indeterminacy set of $\hat f$ as is denoted 
as $\cali_{\hat f}$. $\Gamma_{\hat f}$ is called the graph of the meromorphic 
mapping $\hat f$ or, simply a {\slsf meromorphic graph} if $\hat f$ is either 
irrelevant or clear from the context.

\smallskip By an analytic set in an open subset $\cald$ of a Hilbert
manifold we understand a closed subset $\cala$ of $\cald$ such that 
for every point $a\in\cala$ there exists a neighbourhood $\calv\ni a$ 
and an ideal $\calj$ in the ring $\calo^c_{\calv}$ of bounded holomorphic
functions in $\calv$ such that $\cala\cap\calv = V(\calj)$.

\medskip\noindent{\slsf 1. Finite dimensional case.} When $\calx $ and $\caly$ 
are finite dimensional our definition coincides with that which was given above
and is due to Remmert, see \cite{Re}. Indeed, a subset of a finite dimensional 
connected Riemann manifold is bounded and complete if and only if it is compact.
In that case our assumptions mean that $\pi_{\calx}|_{\Gamma_{\hat f}}: 
\Gamma_{\hat f}\to \calx$ is a proper modification.

\smallskip\noindent{\slsf 2.  Projection $\l^2\setminus\{0\}\to \pp (\l^2)$
revisited.} Notice that $\Gamma_{\widehat{\pr}}$ from Example \ref{hilb-proj} is meromorphic in the sense of our definition with respect to the standard 
Hermitian metric on the base and the Fubini-Study metric on the projective space.
One can easily see this taking into account the fact that $\pp (\l^2)$ is itself 
bounded and complete in the aforementioned metric. 

\smallskip\noindent{\slsf 3. Metrics.} Our definition  is local with 
respect to the base $\calx$ and therefore we can assume that $\calx$ is the unit 
ball in the complex Hilbert space  $L=\cc^n$ or $\l^2$ equipped with the standard 
Hermitian metric. But it, \ie definition, does depend on a metric on the target 
manifold $\caly$. In the case $\caly = \pp (\l^2)$ we always assume the metric 
to be the Fubini-Study one.

\newprg[INT.proj]{Mappings to the projective Hilbert  space}

Now let us formulate the main result of this paper, it concerns with the case 
when $\caly $ is the projective Hilbert  space $\pp (\l^2)$. In \cite{AI} 
the following was proved. Let $\cali$ be an analytic set of codimension 
$\ge 2$  in a complex Hilbert manifold $\calx$ and let $f:\calx\setminus \cali
\to \pp (\l^2)$ be a holomorphic map. Assume that $f(\calx)\not\subset \calh_0$ and set $\calx_0 \deff f^{-1}(\calu_0)$. Write  the holomorphic map 
$f|_{\calx_0}:\calx_0\to \calu_0\equiv \l^2$  as $f(z) = (f_1(z),f_2(z),...)$. 

\begin{nnthm}
\label{lp-mer1}
In these circumstances the following holds:
\begin{enumerate}
\item[1)]  all holomorphic functions $f_k$ extend as meromorphic functions 
to the whole of $\calx$;

\smallskip\item[2)]
for every point $a\in \calx$ there exists a neighbourhood $\calv\ni a$ 
and holomorphic in $\calv$ functions $\phi_0,\phi_1,...$ without common factors 
such that $\sum_{j=0}^{\infty}|\phi_j(z)|$ is locally bounded and such that
on $\calv\setminus I$ one has $ f(z) = [\phi_0(z):\phi_1(z):...]$;

\smallskip\item[3)] in addition
$I\cap \calv \supset \{z\in \calv:\phi_0(z) = \phi_1(z)=...=0\}$;
 
\smallskip\noindent
\item[4)] moreover, the restrictions of meromorphic functions $f_1,f_2,...$  
to $\calv$ are equal to  $\frac{\phi_1}{\phi_0},\frac{\phi_2}{\phi_0},...$ respectively.
\end{enumerate}
\end{nnthm}

\begin{nnrema} \rm
\smallskip\noindent{\bf a)} These statements imply, in particular, that the 
order of poles of $f_k$-s stabilize along $\cala_0\deff f^{-1}(\calh_0)$.
Assumption  $f(\calx)\not\subset \calh_0$ is not essential since 
$\bigcap\calh_k = \emptyset$.

\smallskip\noindent{\bf b)} It should be underlined here that $f|_{\calx_0}:
\calx_0\to \calu_0$ and $\phi \deff (\phi_0,\phi_1(z),...)$ are holomorphic
mappings with values in $\l^2$ and therefore their $\l^2$-norms $\norm{f(z)}^2 
= \sum_{k=1}^{\infty}|f_k(z)|^2$ and $\norm{\phi (z)}^2 = 
\sum_{k=0}^{\infty}|\phi (z)|^2$ are locally bounded.

\smallskip\noindent{\bf c)} As a conclusion we see that the indeterminacy 
set $\cali_{\hat f}$ of our $\hat f$ is locally $\cali_{\hat f}\cap \calv =
\{z\in \calv:\phi_k(z)=0 \,\, \forall k\}$.

\smallskip\noindent{\bf d)} And finally, since a meromorphic {\slsf in any 
sense(!)} mapping should be holomorphic out of codimension $\ge 2$ this
theorem requires that the "extended" by this theorem map should be called
meromorphic!
\end{nnrema}
Therefore we need to prove that the graph of such $f$ as in (1)-(4)
of Theorem \ref{lp-mer1} closes to a meromorphic graph $\Gamma_{\hat f}$
in the sense of the Definition \ref{inf-dim} above.

\begin{nnthm}
\label{lp-mer2}
Let $f$ be as in (1)-(4) of Theorem \ref{lp-mer1}. Then there exists an 
irreducible analytic set $\Gamma_{\hat f}$ in $\calx\times \pp (\l^2)$ 
which satisfies items (1) and (2) of the Definition \ref{inf-dim}, 
\ie $f$ defines a meromorphic mapping $\hat f:\calx \mto \pp (\l^2)$.
Moreover, for $\calv$ and $\phi_k\in \calo (\calv)$ as in item (2) of 
Theorem \ref{lp-mer1} one has  $\cali_{\hat f}\cap \calv  = \{z\in \calv:
\phi_0(z) = \phi_1(z)=...=0\}$. 
\end{nnthm}

Let us notice that in our case $\caly=\pp (\l^2)$ item (1) of the 
Definition \ref{inf-dim} follows from item (2) since, as it was already 
mentioned, $\pp (\l^2)$ is itself complete and bounded. Therefore all
we need to prove is the existence of an 
{\slsf irreducible} analytic set extending the graph. 

\newprg[INT.m-maps]{Irreducibility of the graph}

Now let us try to describe the equation of the meromorphic graph with values 
in $\pp (\l^2)$. We mean the following: if $f$ is locally represented as in 
Theorem \ref{lp-mer1} then as $\Gamma_{\hat f}$ one may naturally take the 
following analytic set 
\begin{equation}
\eqqno(gr-eq1)
\Gamma \deff \{(z,[w])\in \calv\times \pp (\l^2): w_0\phi_k(z) - 
w_k\phi_0(z) = 0 \quad\text{for} \quad k=1,... \infty \}.
\end{equation}
Indeed, we see that $\Gamma_f = \Gamma\cap 
\left[(\calv\setminus I_f)\times \pp (\l^2)\right]$ and $\overline{\Gamma}_f 
\subset \Gamma$. But the latter inclusion might be proper! That's to say 
the meromorphic graph $\Gamma_{\hat f}$ is an irreducible component 
of $\Gamma$ and not the whole $\Gamma$ in general. The  case of the 
canonical projection $\pr : \l^2\setminus \{0\} \to \pp (\l^2)$ considered
above, when actually $\Gamma=\bar\Gamma_f$, is not a typical one even in 
finite dimensions. Let us give an example. 

\begin{nnexmp} \rm
\label{mer-gr1}
Consider the meromorphic map $\hat f:\cc^2 \mto \pp^2$ defined 
on $\cc^2 \setminus \{0\}$  as
\[
f(z) = [\phi_0(z):\phi_1(z):\phi_2(z)] = [z_1+z_2:z_1^2+z_2^2:z_1^3+z_2^3].
\]
We have that $I_{\hat f} = \{0\}$. The analytic set $\Gamma$, defined as 
in \eqqref(gr-eq1), is equal to the graph $\Gamma_f$ of the holomorphic 
part of  $f$ plus $\{0\}\times \pp^2$. But let us see that this is more
than the closure of $\Gamma_f$. The representation of $\hat f$ by 
meromorphic functions is
\[
\left(f_1(z),f_2(z)\right) = \left(\frac{z_1^2+z_2^2}{z_1+z_2}, z_1^2-z_1z_2+z_2^2\right)
\]
and therefore the graph $\Gamma_f$ has the following equation in $(\cc^2\setminus\{0\})\times \pp^2$
\[
\Gamma_f  = \{(z,[w]): w_0(z_1^2+z_2^2) - w_1(z_1+z_2)=0, 
w_2-w_0(z_1^2 -z_1z_2+z_2^2)=0\}.
\]
The same equation defines in $\cc^2\times \pp^2$ the irreducible analytic 
set which is the closure $\bar\Gamma_f$. It is straightforward to check
that $\bar\Gamma_f\cap (\{0\}\times \pp^2) = H_2$, 
where $H_2$ stands for the line $\{w_2=0\}$ in $\pp^2$. 
\end{nnexmp}

\smallskip In finite dimensions this difficulty can be overcome simply 
by taking as $\Gamma_{\hat f}$ the irreducible component of 
$\Gamma$ which contains $\Gamma_f$. But in infinite dimensions the notion 
of an irreducible component is not well defined in general. Therefore our 
goal in this paper, let us say this once more, is to write down an 
{\slsf irreducible} equation for $\Gamma_{\hat f}$ {\slsf explicitly}.
More precisely: we need to write down an equation of an analytic set
$\Gamma_{\hat f}$ in $\calv\times \pp (\l^2)$ such that:
\begin{enumerate}
\smallskip\item[\bf (a)] $\Gamma_{\hat f}$ is  {\slsf irreducible}, that is $\Gamma_{\hat f}$
is an irreducible branch of $\Gamma$ containing $\Gamma_f$\,;

\smallskip\item[\bf (b)] its intersection with $(B\setminus \cali_{\hat f})
\times \pp (\l^2)$ coincides with $\Gamma_f$\,; 

\smallskip\item[\bf (c)] $\Gamma_{\hat f} = \overline{\Gamma}_f$ in classical 
topology. 
\end{enumerate}

\smallskip Item {\bf (a)} will be achieved in a rather explicit way by constructing 
a {\slsf prime} ideal $\widehat{\calp}$ in an appropriate ring of homogeneous polynomials
with holomorphic coefficients, whose zero locus is $\Gamma_{\hat f}$. 
Item {\bf (b)} will be checked  in a rather 
straightforward way. What concerns item {\bf (c)} it will be achieved via
a certain {\slsf finite-codimensional} approximation. 

\newsect[MERM]{Prime ideals and meromorphic graphs}

\newprg[MERM.krull]{Equation of the graph}
 The notation 
$\l^2$  stands for the Hilbert space of sequences of complex numbers 
$z=\{z_k\}_{k=1}^{\infty}$ such that $||z||^2:=\sum_k|z_k|^2<\infty $ 
with the standard Hermitian scalar product $(z,w) = \sum_kz_k\bar w_k$ 
and the standard basis $\{ e_1,e_2,...\}$. 
We are working in a neighbourhood $\calv$ of a point $z_0\in\l^2$ such that
on $\calv\setminus \cali_{\hat f}$ our holomorphic mapping $f$ is defined
by the reduced representation $[\phi_0:\phi_1:... ]$. Here $\phi_k$ are holomorphic in $\calv$ without common factors and $\norm{\phi (z)}^2 = 
\sum_{k=0}^{\infty}|\phi_k(z)|^2$ is locally bounded. The indeterminacy set 
$\cali_{\hat f}$ of the "future" meromorphic map $\hat f$ is defined here
as $\cali_{\hat f}\cap \calv \deff \{z\in \calv : \phi_0(z)=\phi_1(z) = ...=0\}$. 
We may assume that $\phi_0 \not\equiv 0$ and set $\cali_0 = \{z:\phi_0(z) = 0\}$. 

\smallskip Fix a point  $z_0\in \cali_{\hat f}$, \ie such that all $\phi_j$ 
vanish at $z_0$. Notice that our problem is local near every such $z_0$.
Take a point $a\in \calv \setminus \cali_0$ sufficiently close to $z_0$ in order 
to find a ball centred at $a$ which contains $z_0$
and such that its closure is contained in $\calv$. Denote this ball again 
by $\calv$ and by $z_1,z_2,...$ denote the coordinates on $\calv$ centred at 
$a$ such that $\calv$ is the unit ball in these coordinates. Our 
first task is to find such an irreducible analytic set $\Gamma_{\hat f}$
in $\calv\times \pp (\l^2)$ that its intersection with $(\calv\setminus 
\cali_0)\times \pp (\l^2)$ coincides with $\Gamma_f\cap \left[(\calv
\setminus \cali_0)\times \pp (\l^2)\right]$. Notice that then, obviously 
$\Gamma_{\hat f}\cap (\calv \setminus \cali_{\hat f})\times \pp (\l^2)$ 
will coincide with $\Gamma_f$.

\smallskip 
Restrict equation \eqqref(gr-eq1) to $(\calv\setminus \cali_0)\times \calu_0$ to obtain 
the equation of the graph of the holomorphic mapping $f|_{\calv\setminus \cali_0}:
\calv\setminus \cali_0 \to \calu_0$, which in affine coordinates $\frac{w_k}{w_0}$ 
on $\calu_0$ writes as 
\begin{equation}
\eqqno(gr-eq3)
\begin{cases}
\frac{w_1}{w_0} \phi_0(z) - \phi_1(z) = 0, \cr
\frac{w_2}{w_0} \phi_0(z) - \phi_2(z) = 0, \cr
...............................
\end{cases}
\end{equation}

Denote these coordinates as $\vect_j\deff \frac{w_j}{w_0}$ and set

\begin{equation}
\eqqno(def-of-g)
\begin{cases}
g_1(z)[\vect] \deff \vect_1 \phi_0(z) - \phi_1(z),\cr
g_2(z)[\vect] \deff \vect_2 \phi_0(z) - \phi_2(z),\cr
.........................................
\end{cases} 
\end{equation}

 We obtain the equation of the affine part $\Gamma_0$ of $\Gamma$, \ie of 
$\Gamma \cap [(\calv\setminus \cali_0)\times \calu_0]$ as follows

\begin{equation}
\Gamma_0 = \{(z,\vect)\in (\calv\setminus \cali_0)\times \calu_0: g_1(z)[\vect_1]=...
=g_k(z)[\vect_k]= ...= 0\}.
\end{equation}

\medskip
This means that  $\Gamma_0=\Gamma_f\cap \left[(\calv\setminus \cali_0) \times 
\calu_0)\right] = V(g_1,g_2,...)\cap \left((\calv\setminus \cali_0)\times 
\calu_0\right)$. Here we denote by  $V(g_1,g_2,...)$ 
the zero locus of the ideal $(g_1,g_2,...)$ generated by $g_1,g_2,...$ 
in the ring $\calo^c_{\calv}[\vect]$ of continuous polynomials in $\vect = 
(\vect_1,\vect_2,...)$ with coefficients in $\calo^c_{\calv}$. The latter 
stands for the ring of germs of  bounded holomorphic functions 
in $\calv$.  Notice that $g_k$-s are such because such are $\phi_k$-s.

\newprg[MERM.prime]{Prime ideals in non-Noetherian rings}
Consider the following ideal in $\calo^c_{\bar \calv}[\vect]$:
\begin{equation}
\eqqno(ideal-1)
\calp = \left\{g\in \calo^c_{\calv}[\vect]: g = \sum \frac{p_k}{q}g_k, \text{ where } q, p_k \in 
\calo^c_{\calv}[\vect] \text{ with } q(0)[0]\not=0\right\},
\end{equation}
where sums are finite. 
\begin{thm}
\label{prime-id1}
Ideal $\calp$ is prime and $V(\calp)\cap [( \calv\setminus \cali_0)\times 
\calu_0] = \Gamma_f\cap [( \calv\setminus \cali_0)\times \calu_0]$. 
\end{thm}
\proof Denote by $\calo^c_{\calv}[\vect]_{(z,\vect)}$ the localization of 
the ring $\calo^c_{\calv}[\vect]$ by the multiplicatively closed set 
$\calo^c_{\calv}[\vect]\setminus (z,\vect)$, here $(z,\vect)\deff (z_1,...,
z_k,...\vect_1,...,\vect_k,...)$ denotes the ideal generated by $z_1,...,z_k,
...\vect_1,...,\vect_k,...$. Consider the following ideal in $\calo^c_{\calv}
[\vect]_{(z,\vect)}$:
\begin{equation}
\eqqno(ideal-2)
\calp' = \left\{g\in \calo^c_{\calv} [\vect]_{(z,\vect)}: g = 
\sum\frac{p_k}{q}g_k, \text{ where } q, p_k \in 
\calo^c_{\calv} [\vect] \text{ with } q(0)[0]\not=0\right\}.
\end{equation}
In other words $\calp' = (g_1,g_2,...)$ is the ideal generated by 
$g_1,g_2,...$ in the ring $\calo^c_{\calv} [\vect]_{(z,\vect)}$. It is
clear that $\calp'\cap \calo^c_{\calv} [\vect] =\calp$ and therefore 
$\calp$ will be prime if we prove that $\calp'$ is prime. Consider the 
ring $\cc \{z\}[[\vect]]$ of formal power series in  $\vect$ with coefficients 
in the ring $\cc \{z\}$ of convergent power series in $z$. Furthermore, 
consider the following ideal in $\cc \{z\}[[\vect]]$:
\begin{equation}
\eqqno(ideal-3)
\calp^{''} = \left\{g\in \cc \{z\}[[\vect]]: g = \sum p_kg_k, \text{ where } 
p_k \in \cc \{z\}[[\vect]]\right\}.
\end{equation}
I.e., $\calp^{''} = (g_1,g_2,...)$ is again the ideal generated by 
$g_1,g_2,...$, but now in $\cc \{z\}[[\vect]]$. We have the natural inclusion
$\calo^c_{\calv}[\vect]_{(z,\vect)}\subset \cc \{z\}[[\vect]]$ given by 
the power series expansions of the elements of $\calo^c_{\calv}[\vect]_{(z,\vect)}$ at zero:
\begin{equation}
\eqqno(pow-exp1)
\left(1-\sum_{p=0}^d\sum_{|J|=p}q_J(z)\vect^J\right)^{-1} = 1 + \sum_{n=1}^{\infty}\left(\sum_{p=0}^d\sum_{|J|=p}
q_J(z)\vect^J\right)^n = 1 + \sum_KQ_K(z)\vect^K,
\end{equation}
where
\begin{equation}
\eqqno(pow-exp2)
Q_K(z) = \sum_{J_1+...+J_l=K}q_{J_1}(z)...q_{J_l}(z).
\end{equation}
I.e., all $Q_K$ are finite sums  and therefore the  series in the right 
hand side of \eqqref(pow-exp1) belong to $\cc \{z\}[[\vect]]$. Here $K$ 
runs over the set of all multi indexes of finite length, \ie $K\in 
\nn_0^{\infty}\deff\bigcup_{k=1}^{\infty} \nn^k\cup \{0\}$. Notice 
that since $q_0(0) =0$ in \eqqref(pow-exp1) so $Q_0(0)=0$ as well.

\begin{lem} 
\label{prime-id2}
$\calp^{''}\cap \calo^c_{\bar B} [\vect]_{(z,\vect)} = \calp'$.
\end{lem}
\proof Inclusion $\calp' \subset \calp^{''}\cap \calo^c_{\bar B} [\vect]_{(z,\vect)}$ is obvious. To prove the opposite, one notices that 
\[
\calp^{''}\cap \calo^c_{\bar B} [\vect]_{(z,\vect)} = \Big\{g\in \calo^c_{\bar B}[\vect]_{(z,\vect)}: 
g=\sum p_kg_k, p_k\in \cc\{z\}[[\vect]]\Big\} \subset \bigcap_{N=1}^{\infty}\Big\{g\in 
\calo^c_{\bar B}[\vect]_{(z,\vect)}:
\]
\[
g = \sum p_kg_k +r\text{ with } p_k\in \cc \{z\}[\vect], \text{ } \deg_{\vect} p_k \le N-1,
\text{ and remainder } r \text{ vanishing up }
\]
\begin{equation}
\eqqno(inclus1)
\text{   to  terms of order } \ge N \text{ in } \vect \Big\} = \bigcap_{N=1}^{\infty} \left\{\calp' + 
(\vect)^N\cap \calo^c_{\bar B}[\vect]_{(z,\vect)}\right\} = \calp'. \qquad\qquad\qquad
\end{equation}
As for the inclusion above notice that if $g\in \calo^c_{\bar B}[\vect]_{(z,
\vect)}\cap \calp^{''}$ then $g$ writes as $g = \sum p_kg_k$ for some $p_k\in 
\cc \{z\}[[\vect]]$. Monomials of order $\ge N$ in $p_k$ give contribution 
to $r$ and we can forget about them. I.e., we can assume that $\deg_{\vect} 
p_k\le N-1$. The last equality follows from the obvious fact that $\bigcap_{N=1}^{\infty}(\vect)^N=(0)$ in 
$\cc\{z\}[[\vect]]$ and therefore $\bigcap_{N=1}^{\infty}(\vect)^N\cap \calo^c_{\bar B}[\vect]_{(z,\vect)}=(0)$
in $\calo^c_{\bar B}[\vect]_{(z,\vect)}$ viewed as a subring of $\cc\{z\}[[\vect]]$ as above. The argument 
is standard: in the factor ring $\calo^c_{\bar B}[\vect]_{(z,\vect)}/\calp'$ the image 
$\calj\deff \left[(\vect)\cap \calo^c_{\bar B}[\vect]_{(z,\vect)}\right]/\calp'$ satisfies $\bigcap_{N=1}^{\infty}
\calj^N = (0)$ and therefore $\bigcap_{N=1}^{\infty} \left\{\calp' + 
(\vect)^N\cap \calo^c_{\bar B}[\vect]_{(z,\vect)}\right\}$ is nothing but $\calp'$.


\smallskip 

\smallskip\qed

\smallskip  Therefore to prove that $\calp'$ is prime it is sufficient to prove that $\calp^{''}$ is prime.

\begin{lem} 
$\calp^{''}$ is prime.
\end{lem}
\proof  Recall that $\calp^{''} = (g_1,...,g_k,...)$ in $\cc\{z\}[[\vect]]$. All we need to prove 
is that $\cc\{z\}[[\vect]]/(g_1,...,g_k,...)$ is an integral domain. Let's start from $\cc\{z\}[[\vect]]/(g_1)$. 
Since $g_1(z)[\vect]$ is a Weierstrass polynomial of degree one with respect to $\vect_1$ we get by 
the Weierstrass division theorem  that any formal power series $f(z,\vect)$ uniquely writes as 
\begin{equation}
\eqqno(weier-1)
f(z,\vect) = a(z,\vect)g_1(z,\vect) + b(\vect_2,...,\vect_k,...).
\end{equation}
This means that $\cc\{z\}[[\vect]]/(g_1)$ is isomorphic to $\cc\{z\}[[\vect_2,...,\vect_k,...]]$. 
By induction, using the fact that every $g_j$ is a Weierstrass polynomial w.r.t. $\vect_j$ we prove that
the ring $\cc\{z\}[[v]]/(g_1,...,g_N)$ is isomorphic to $\cc\{z\}[[\vect_{N+1}, ...]]$ for every $N$. 
And this means that $\cc\{z\}[[\vect]]/(g_1,...,g_k,..)$ is isomorphic to $\cc\{z\}$. The latter is an 
integral domain. Therefore $\calp^{''} = (g_1,...,g_k,...)$ is prime and our 
lemma together with the first assertion of the theorem are proved.

\smallskip  To prove the second assertion of the theorem consider the 
analytic set $V(\calp)$ in $ \calv\times \calu_0$. We need to show that 
\begin{equation}
\eqqno(gamma=vp)
\Gamma_f\cap \left[( \calv\setminus \cali_0)\times \calu_0\right] = V(\calp)\cap 
\left[( \calv\setminus \cali_0)\times \calu_0\right].
\end{equation}
Since $\calp \ni g_1,...,g_k,...$ we have that $V(\calp)\subset 
V(g_1,...,g_k,...)$ and the intersection of the latter with $( \calv\setminus 
\cali_0)\times \calu_0$ is clearly $\Gamma_f\cap [( \calv\setminus \cali_0) \times 
\calu_0]$. Furthermore, notice that germs at $(0,0)\in  \calv\times \calu_0$
of analytic sets $V(\calp)$ and $\Gamma_f$ coincide. This is because the
factors $q$ in the definition of $\calp$ are required not to vanish at the 
origin. Therefore, since $\Gamma_f$ is a manifold in $(\calv\setminus \cali_0)\times 
\calu_0$, one must have that $V(\calp)\cap \left[(\calv\setminus \cali_0)\times 
\calu_0\right] = \Gamma_f\cap \left[(\calv\setminus \cali_0)\times \calu_0\right]$. 
Theorem \ref{prime-id1} is proved.

\smallskip\qed

With an ideal $\calp$ in $\calo^c[\vect]$ we can associate a {\slsf homogeneous}
ideal 
\begin{equation}
\eqqno(hat-id1)
\widehat{\calp} \deff \left\{ \hat f(\wect_0,...,\wect_k,...) = 
\wect_0^{\deg f}f\left(\frac{\wect_1}{\wect_0},...,
\frac{\wect_k}{\wect_0},...\right):f\in \calp\right\}
\end{equation}
in $\calo^c_{\calv}[w]$, here $w = \{\wect_0,...,\wect_k,... \}$, see ex.
\cite{Mu}. 
\begin{lem}
\label{prime-id3}
Let $\calp$ be a prime ideal in $\calo^c_{\bar B}[\vect]$. Then  the 
homogeneous ideal $\widehat{\calp}$ is prime. Moreover, its zero locus 
$V(\widehat{\calp})$ coincides with the closure
of $V(\calp)$  in $\calv\times \pp (\l^2)$ .
\end{lem}
\proof If $\hat f\hat g\in \widehat{\calp}$ then this means that 
$\wect_0^{\deg f}f\left(\frac{\wect_1}{\wect_0},...,
\frac{\wect_k}{\wect_0},...\right)\wect_0^{\deg g}g\left(\frac{\wect_1}
{\wect_0},..., \frac{\wect_k}{\wect_0},...\right) = $ 
$=\wect_0^{\deg f+\deg g}(fg)\left(\frac{\wect_1}{\wect_0},...,
\frac{\wect_k}{\wect_0},...\right)$ with $fg\in \calp$. The primality 
of $\widehat{\calp}$ follows from that of $\calp$. Since $\calu_0=\{[\wect]:\wect_0\not=0\}$ the intersection of the zero set of every $\hat f$ 
from $\widehat{\calp}$ with $\calv\times \calu_0$ coincides with the zero set of 
$f$, here $f$ and $\hat f$ are related as in \eqqref(hat-id1). And the former
is the closure of the latter, see Lemma II.4.1.1 on the page 80 in \cite{Ra}.
Now 
\begin{equation}
\eqqno(inters-1)
V(\widehat{\calp}) = \bigcap_{f\in \calp} V(\hat f) = \bigcap_{f\in \calp} 
\overline{V(f)} = \overline{\bigcap_{f\in \calp} V(f)} =\overline{V(\calp)}.
\end{equation}
Lemma is proved.

\smallskip\qed

\smallskip Items {\bf (a)} and {\bf (b)} being proved we turn
to item {\bf (c)}.




\newsect[APPRX]{Finite codimensional approximation}
\newprg[APPRX.graph]{Graphs of infinite codimension}
Analytic set $V(\widehat{\calp})$ which pretends to be the closure of the 
graph of our mapping $f$, \ie the graph $\Gamma_{\hat f}$ of the resulting 
meromorphic map is of infinite both dimension and codimension.
Such analytic sets look very differently from the finite dimensional
analytic sets. Consider the following simple examples.

\begin{exmp} \rm
Let $G_1$ denote the Banach space $l^1$ of complex sequences
$z=\{z_n\}_{n\ge 1}$ with finite $l^1$-norm $\norm{z}_1 = \sum_n|z_n|$. 
Further denote by $G_2$ the Banach space of double sequences $\{z_{jk}\}_{1\le j<k}$ with the norm $\norm{z}_2 = \sum_{jk}|z_{jk}|$. Set
$G\deff G_1\oplus G_2$ and take as the norm the sum of corresponding
norms. Consider the following holomorphic mapping $f:l^2\to G$:
\begin{equation}
f : \{z_n\} \to \Big(\{z_n(z_n - 1/n)\}, \{z_jz_k\}\Big).
\end{equation}
The zero set of $f$ is $\scra = \{0, e_1, \frac{1}{2}e_2,..., 
\frac{1}{n}e_n, ...\}$, where $<e_1,...,e_n,...>$ is the standard
orthonormal basis of $l^2$. Therefore this analytic set is, in 
fact, a sequence of points in $l^2$ converging to zero.
\end{exmp}

The situation can be even worse. 

\begin{exmp}\rm
Consider the following meromorphic mapping $\hat f:\l^2\mto \pp (\l^1)$
\begin{equation}
\eqqno(indet1)
\hat f : z=\{z_k\}_{k=1}^{\infty} \mapsto \left[\{z_k(z_k-1/k)\}_{k=1}^{\infty}\right],
\end{equation}
or, in more details 
\[
\{z_k\}_{k=1}^{\infty} \mapsto  \left[z_1(z_1-1):...: z_k(z_k-1/k):...
\right].
\]	
The indeterminacy set $\cali$ of this $\hat f$ (\ie all $z$ such that all 
components in \eqqref(indet1) vanish) is of Cantor type, see pp. 33-34 
in \cite{Ra}. In particular it is totally disconnected and (again)
contains the sequence $\{e_1,...,\frac{1}{k}e_k,...\}$ converging to 
zero.
\end{exmp}

\smallskip These examples are not the ones we expect from a
meromorphic graph.

\newprg[APPRX.fin-def]{Graphs of finite definition}

However analytic sets of finite codimension (they are mostly called
of {\slsf finite definition}) in Hilbert manifolds possess the properties
analogous of that ones of finite dimension. This was proved by Ramis,
and we shall consequently use results from his book \cite{Ra}. We 
say that an analytic set $\cala$ in an open subset $\cald$ of 
a Hilbert manifold is of {\slsf finite definition} if ideals, 
which define $\cala$ locally, are finitely generated.

\smallskip Let $g_k(z)[\vect]\in \calo^c_B[\vect], k=1,2,...$ be 
polynomials as in \eqqref(def-of-g). Analogously to \eqqref(ideal-1)
set 
\begin{equation}
\eqqno(ideal-n1)
\calp_n \deff  \left\{g\in \calo^c_{\bar B}[\vect]: g = 
\sum_{k=1}^n \frac{p_k}{q}g_k, \text{ where } q, p_k \in 
\calo^c_{\bar B}[\vect] \text{ with } q(0)[0]\not=0\right\}.
\end{equation}

We observe the following items.

\smallskip
\begin{itemize}
\smallskip\item $\calp_n$ and, consequently, the corresponding homogeneous
ideal $\widehat{\calp}_n$ are prime.

\item Analytic sets $V\left(\widehat{\calp}_n\right)\cap 
\big[(B\setminus \cali_0)\times\pp (\l^2)\big]$ are of finite definition.

\item $V(\widehat{\calp}_n)$ is the closure of  $V\left(\widehat{\calp}_n\right)
\cap \big[(\calv\setminus \cali_0)\times\pp (\l^2)\big]$.

\smallskip\item Moreover, one has $\widehat{\calp}_n \subset 
\widehat{\calp}_{n+1}$ and $\bigcup \widehat{\calp}_n = \widehat{\calp}$.
And therefore 
\begin{equation}
\eqqno(inters-2)
V\left(\widehat{\calp}_n\right)\supset V\left(\widehat{\calp}
_{n+1}\right) \quad\text{ and } \quad \bigcap V(\widehat{\calp}_n) = 
V(\widehat{\calp}).
\end{equation}
\end{itemize}

The first item can be proved by repeating literally the proof of 
Theorem \ref{prime-id1} and Lemma \ref{prime-id3} above. As for
the second item set
\begin{equation}
\eqqno(gamma-n)
\Gamma_n = \{(z;[w]) \in (B\setminus \cali_0)\times \pp (\l^2)): 
\hat g_1(z)([w])=...=\hat g_n(z)([w]) = 0\},
\end{equation}
here $\hat g_k(z)([w])$ are projectivizations of $g_k(z)[\vect]$
as in \eqqref(hat-id1). One has that $\Gamma_n = V\left(\widehat{\calp}_n\right)\cap \big[(B\setminus \cali_0)\times\pp (\l^2)\big]$ for 
the same reason as in \eqqref(gamma=vp), which proves the assertion.
By Theorem II.4.1.2 of \cite{Ra} the closure $\overline{\Gamma}_n$ is 
again of finite codimension in $\calv\times \pp (\l^2)$. Knowing that
$V(\widehat{\calp}_n)$ contains $\overline{\Gamma}_n$ suppose that 
$V(\widehat{\calp}_n)\not\subset\overline{\Gamma}_n$. But then $V(\widehat{\calp}_n)$ is the union of two distinct components $\overline{\Gamma}_n$ and 
$V(\widehat{\calp}_n)\cap \cali_0\times \pp (\l^2)$, contradiction. 
I.e., $V(\widehat{\calp}_n)$ is the closure of  $V\left(\widehat{\calp}_n\right)\cap \big[(\calv\setminus \cali_0)\times\pp (\l^2)\big]$. This gives the second item.

\smallskip The last item ( is obvious and obviously implies that $V(\widehat{\calp})$ is the closure of $\bigcap V\left(\widehat{\calp}_n\right)\cap \big[(B\setminus I_0)\times\pp (\l^2)\big]$. But the 
latter set is $\bigcap\Gamma_n = \Gamma_f$. I.e., $V(\widehat{\calp})
= \overline{\Gamma}_f$. Item {\bf (c)} is proved, \ie our meromorphic 
graph $\Gamma_{\hat f}$ is actually $V(\widehat{\calp})$.

\smallskip

\ifx\undefined\bysame
\newcommand{\bysame}{\leavevmode\hbox to3em{\hrulefill}\,}
\fi

\def\entry#1#2#3#4\par{\bibitem[#1]{#1}
{\textsc{#2 }}{\sl{#3} }#4\par\vskip2pt}

\end{document}